%

%
\magnification=\magstep1
\input amstex
\documentstyle{amsppt}

\def\c{\cite}
\def\range{\text{range}}

\topmatter
\title
On Countably Closed Complete Boolean Algebras
\endtitle
\author
Thomas Jech and Saharon Shelah
\endauthor
\affil
The Pennsylvania State University  \\
The Hebrew University and Rutgers University
\endaffil
\address 
\endgraf Department of Mathematics, The Pennsylvania State University,
University Park, PA 16803, USA 
\endgraf Institute of Mathematics, The Hebrew University, Jerusalem, Israel, and  
Department of Mathematics, Rutgers University, New Brunswick, NJ 08903, USA
\endaddress
\email 
jech\@math.psu.edu, shelah\@sunrise.huji.ac.il, shelah\@math.rutgers.edu 
\endemail
\thanks
The first author has been partially supported by the
U.S.-Czechoslovakia cooperative grant INT-9016754 from the NSF.\endgraf 
The second author has been partially supported by
the U.S.-Israel Binational Science Foundation.  Publication number 
565\endthanks

\abstract
It is unprovable that every complete subalgebra of a countably
closed complete Boolean algebra is countably closed.
\endabstract

\endtopmatter

\document
\baselineskip 20pt

\subhead
{Introduction}
\endsubhead
A partially ordered set $(P,<)$ is $\sigma$-{\it closed} if every
countable chain in $P$ has a lower bound.  A complete Boolean algebra
$B$ is {\it countably closed} if $(B^+,<)$ has a dense subset that
is $\sigma$-closed.  In \c2 the first author introduced a weaker
condition for Boolean algebras, {\it game-closed}: the second player
has a winning strategy in the infinite game where the two players 
play an infinite descending chain of nonzero elements, and the second
player wins if the chain has a lower bound.  In \c1, Foreman proved
that when $B$ has a dense subset of size 
$\aleph_1$ and is game-closed then $B$ is countably closed.
(By Vojt\'a\v s \c5 and Veli\v ckovi\'c \c4 this holds for every $B$
that has a dense subset of size $2^{\aleph_0}$.)
We show that, in general, it is unprovable that game-closed implies
countably closed.  We construct a model in which a $B$ exists that 
is game-closed but not countably closed.  It remains open whether a
counterexample exists in ZFC.

Being game-closed is a hereditary property:  If $A$ is a complete
subalgebra of a game-closed complete Boolean algebra $B$ then $A$ is
game-closed.  It is observed in \c3 that every game-closed algebra is
embedded in a countably closed algebra; in fact, for a forcing notion
$(P,<)$, being game-closed is equivalent to the existence of a
$\sigma$-closed forcing $Q$ such that $P\times Q$ has a dense
$\sigma$-closed subset.  Hence the statement ``every game-closed
complete Boolean algebra is countably closed'' is equivalent to the
statement ``every complete subalgebra of a countably closed complete
Boolean algebra is countably closed''.

Below we construct (by forcing) a model of ZFC+GCH and in it a
partial ordering $P$ of size $\aleph_2$ such that $B(P)$, the completion of
$P$, is not countably closed, but $B(P\times Col)$ is, where $Col$ is
the L\'evy collapse of $\aleph_2$ to $\aleph_1$ (with countable conditions).

\proclaim
{Theorem} It is consistent that there exists a partial ordering $(P,<)$
such that $B(P)$ is not countably closed but $B(P\times Col)$ is
countably closed.
\endproclaim

\subhead
{Forcing Conditions}
\endsubhead

We assume that the ground model satisfies $GCH$.

We want to construct, by forcing, a partially ordered set $(P,<_P)$ of
size $\aleph_2$ that has the desired properties.  We shall use as forcing
conditions countable approximations of $P$.  One part of a forcing
condition will thus be a countable partial ordering $(A,<_A)$ with the
intention that $A$ be a subset of $P$ and that the relation $<_A$ on $A$ be
the restriction of $<_P$.  As $P$ will have size $\aleph_2$, we let
$P=\omega_2$, and so $A$ is a countable subset of $\omega_2$.

The second part of a forcing condition will be a countable set
$B\subset A\times Col$, a countable approximation of a dense set in
the product ordering $P\times Col$.  The third part of a forcing
condition will be a countable set $C$ of countable descending chains
in $A$ that have no lower bound. Finally, a forcing condition includes
a function that guarantees that the limit of the $B$'s is $\sigma$-closed 
(and so $P\times Col$ has a $\sigma$-closed dense subset).

Whenever we use $<$ without a subscript, we mean the natural ordering
of ordinal numbers.

\proclaim
{Definition} \ {\rm For any set $X,$ $Col(X)\,$ is the set of all countable 
functions $q$ such that $dom (q)\in \omega_1$ and range $(q)\subset 
X;$ $Col = Col(\omega_2)$.} 
\endproclaim

\bigpagebreak

\proclaim
{Definition} \ {\rm The set $R$ of forcing conditions $r$ consists of
quadruples $r=((A_r,<_r)$, $B_r$, $C_r$, $F_r)$ such that
\roster
\item $A_r$ is a countable subset of $\omega_2$,
\item $(A_r,<_r)$ is a partially ordered set,
\item if $b<_r a$ then $a<b$,
\item $B_r$ is a countable subset of $A_r\times Col(A_r)$, and for every
	$(p,q)\in B_r$,\newline 
        $p\in \range(q)$,
\item $C_r$ is a countable set of countable sequences $\{a_n\}_{n=0}^{\infty}$
	in $A_r$ with the property that $a_0>_r a_1>_r\cdots >_r a_n>_r\cdots$
	and that $\{a_n\}_n$ has no lower bound in $A_r$,
\item $F_r$ is a function of two variables, $\{a_n\}_n \in C_r$ and 
        $(p,q)\in B_r$ such that $p \ge a_0,$ and $\range(F_r)\subset\omega.$
        If $m=F_r(\{a_n\}_n,(p,q))$ then for every $(p',q')\in B_r$ stronger
        than $(p,q)$, 
$$\text{if }p'<_r a_m \text{ then } p'\perp_r \{a_n\}_n \text{ (i.e.
        $p'\perp_r a_k$ for some $k$).}\tag*$$ 

	If $r,s\in R$ then $r<_R s$ ($r$ is stronger than $s$) if
\item $A_r\supseteq A_s$,
\item $<_r$ and $<_s$ agree on $A_s$, and $\perp_r$ and $\perp_s$ agree on $A_s$; 
        i.e. if $a,b\in A_s$ then $a<_r b$ iff $a<_s b$ and
        $a\perp_r b$ iff $a\perp_s b$ for all $a,b\in A_s]$,
\item $B_r\supseteq B_s$,
\item $C_r\supseteq C_s$,
\item $F_r\supseteq F_s$.
\endroster}
\endproclaim

The relation $<_A$ on $R$ is a partial ordering.  We shall prove that
the forcing extension by $R$ contains a desired example $(P,<_P)$.  That
$R$ is a cardinal-preserving model of $ZFC+GCH$ follows from the next
two lemmas:

\proclaim
{Lemma 1}  $R$ is $\sigma$-closed.
\endproclaim

\demo
{Proof}  Let $\{r_n\}_n$ be a sequence of conditions such that
$r_0>_Rr_1>_R\cdots >_Rr_n>_R\cdots$.  We show that $\{r_n\}_n$ has a
lower bound.

Assuming that for each $n$, $r_n=((A_n,<_n),B_n,C_n,F_n)$, we let
$A_r=\bigcup^\infty_{n=0}A_n$, $B_r=\bigcup^\infty_{n=0}B_n$, 
$C_r=\bigcup^\infty_{n=0}C_n$, $F_r=\bigcup^\infty_{n=0}F_n$
and $<_r=\bigcup^\infty_{n=0}<_n$; we claim that
$r=((A_r,<_r),B_r,C_r,F_r)$ is a condition, and is stronger than each $r_n$.

The triple $r$ has clearly properties (1)--(4). It is also easy to see
that for every $n,$ $<_r$ agrees with $<_n$ and $\perp_r$ agrees with
$\perp_n$ on $A_n.$ To verify (5),
let $\{a_n\}_n\in C_r$.  There is an $m$ such that $\{a_n\}_n\in C_k$
for all $k\geq m$, and therefore $\{a_n\}_n$ has no lower bound in
any $A_k$.  Thus $\{a_n\}_n$ has no lower bound in $A_r$.
Finally, to verify (6), let $F_r(\vec a,(p,q))=m$ and let $(p',q')$ be stronger
than $(p,q).$ Since (*) holds in $r_n$ where $n$ is large enough so that
$\vec a \in C_n$ and $(p,q), (p',q') \in B_n,$ (*) holds in $r$ as well.

Therefore $r$ is a condition and for every $n,$ $r$ is stronger than $r_n.$

\enddemo

\proclaim
{Lemma 2} \ $R$ has the $\aleph_2$-chain condition.
\endproclaim

\demo
{Proof} If $W$ is a set of conditions of size $\aleph_2$, then a
$\Delta$-system argument (using CH) yields two conditions
$r,s\in W$ such that if $r=((A_r,<_r),B_r,C_r, F_r)$ and
$s=~((A_s,<_s),B_s,C_s,F_s),$ then there is a $D$ (the root of the
$\Delta$-system) such that $D=A_r\cap A_s$, $\sup D<\min(A_r-D)$,
$\sup A_r<\min(A_s-D)$, $<_r$ and $<_s$ agree on $D$, 
$\perp_r$ and $\perp_s$ agree on $D$, $B_r\cap(D\times~
Col(D))=B_s\cap(D\times Col(D))$, $C_r\cap D^\omega=C_s\cap D^\omega$,
and $F_r(\vec a,(p,q))=F_s(\vec a,(p,q))$ whenever $\vec a\in C_r\cap
D^\omega$ and $(p,q)\in B_r\cap(D\times Col(D)).$

Moreover, there exists a mapping $\pi$ of $A_s$ onto $A_r$ that is an
isomorphism between $s$ and $r$ and is the identity on $D.$

Let $t=((A_t,<_t), B_t,C_t,F_t)$ where $A_t=A_r\cup A_s$, $B_t=B_r\cup
B_s$, $C_t=C_r\cup C_s$, $<_t=<_r\cup <_s$, and $F_t$ will be defined
below such that $F_t\supseteq F_r\cup F_s$.
We claim that $t$ is a condition, and is stronger than both $r$ and $s$; 
thus $r$ and $s$ are compatible.  Properties (1)--(4) are easy to verify.  
It is also easy to see that $<_t$ agrees with $<_r$ on $A_r$ and with
$<_s$ on $A_s,$ and $\perp_t$ agrees with $\perp_r$ on $A_r$ and with
$\perp_s$ on $A_s.$

Note that if $a\in A_r-D$ and $b\in A_s-D$ then $a\perp_t b.$ Thus if
$\{a_n\}_n$ is in $C_r$ but not in $C_s$ (or vice versa) then
$\{a_n\}_n$ has no lower bound in $A_r\cup A_s$, and so (5) holds.

In order to deal with (6), we first verify it for the values of $F_t$
inherited from either $r$ or $s.$ Thus let $\vec a \in C_r,$ $(p,q)\in B_r$,
$m=F_r(\vec a,(p,q))$ and let $(p',q')\in B_t$ be stronger than $(p,q).$
(The argument for $s$ in place of $r$ is completely analogous.) If
$(p',q')\in B_r$ then (*) holds in $r$ and therefore in $t.$ Thus assume
that $(p',q')\in B_s.$

Since $p'\in A_s$ and $p'<_t p,$ it follows that $p\in D,$ and since $\range(q)
\subseteq \range(q') \subseteq A_s,$ we have $(p,q)\in B_s.$ Now if $\vec a \in
C_s$ then $F_s(\vec a,(p,q))=F_r(\vec a,(p,q))$ and so $p'$ satisfies (*) 
in $s$ and hence in $t.$ If $\vec a \notin C_s$ and $p'\notin A_r$ then
$p' \perp_t \vec a$ and again $p'$ satisfies (*).

The remaining case is when $p'\in D$ and $(p,q)\in B_r \cap B_s.$ Since
$(p',\pi q')=(\pi p',\pi q')$ is stronger than $(p,q)=(\pi p,\pi q),$ $p'$
satisfies (*) in $r$ and therefore in $t.$

To complete the verification of (6) we define $F_t(\vec a,(p,q))$
for those $\vec a$ and $(p,q)$ that come from the two different conditions.
Let $\vec a \in C_r - C_s$ and $(p,q)\in B_s - B_r$ (the other case being
analogous) be such that $p \ge a_0.$ We let $F_t(\vec a,(p,q))$ be the least
$m$ such that $a_m \notin D.$ 

Let $(p',q')\in B_t$ be stronger than $(p,q)$; we'll show that $p'\nless_t a_m.$
This is clear if $p'\in D.$ If $p'\notin D$ then $p'$ cannot be
in $A_r$ because by (4) $p'\in \range(q')\subseteq\range(q)\subseteq A_s.$ 
It follows that $p' \perp_t a_m.$

Therefore $t$ is a condition and is stronger than both $r$ and $s.$

\bigpagebreak

Let $G$ be a generic filter on $R$.  In $V_G$, we let $P=\bigcup
\{A_r:r\in G\}$, $<_P=\bigcup\{<_r:r\in G\}$, and $Q=\bigcup\{B_r:r\in G\}$.
$(P,<_P)$ is a partial ordering and $Q\subset P\times Col$.  We shall
prove that $Q$ is $\sigma$-closed and is dense in $P\times Col$, and
that the complete Boolean algebra $B(P)$ does not have a dense 
$\sigma$-closed subset.
\enddemo

\proclaim
{Lemma 3}  $P=\omega_2$.
\endproclaim

\demo
{Proof} We prove that for every $s$ and every $p\in \omega_2$ there
exists an $r<_R s$ such that $p\in A_r$.  But this is
straightforward: let $A_r=A_s\cup\{p\}$, $B_r=B_s$, $C_r=C_s,$ $F_r=F_s$ and
$<_r=<_s$; properties (1)--(11) are easily verified. (Note that $p \perp_r a$
for all $a\in A_s.)$

\enddemo

\proclaim
{Lemma 4}  $Q$ is dense in $P\times Col$.
\endproclaim

\demo
{Proof}  Let $s$ be a condition and let $p_0\in A_s$ and $q_0\in Col$.
We shall find an $r<_R s$, $p\in A_r$ and $q\supset q_0$ such that $p<_r p_0$ 
and $(p,q)\in B_r$:  Let $p$ be an ordinal greater than $p_0$ and 
$p\notin A_s$, let $q\in Col$ be such that $q\supset q_0$ and $p\in$
range $(q)$, and let $A_r=A_s\cup$ range $(q)$, $B_r=B_s\cup
\{(p,q)\}$, $C_r=C_s$, and let $<_r$ be the partial order of $A_r$
that extends $<_s$ by making $p<_r p_0$.  
Finally, let $F_r(\vec a,(p,q))=0$ for all $\vec a \in C_r.$

To see that $r=((A_r,<_r),B_r,C_r,F_r)$ is a condition, 
note that for every $\vec a \in C_r,$ $p$ is not a lower bound of $\vec a$
(because $p_0$ isn't) and hence $p\perp_r \vec a.$ This implies both (5)
and (6). Since adding $p$ does not effect the relation $\perp$ on $A_s,$
we have (8) and so $r$ is stronger than $s.$ 
\enddemo

Next we prove that $Q$ is $\sigma$-closed.

\proclaim
{Lemma 5}  If $u=\{(p_n,q_n)\}_{n=0}^{\infty}$ is a descending chain in
$Q$ then $u$ has a lower bound.
\endproclaim

\demo
{Proof} Let $\dot{u}$ be a name for a descending chain and let $s$ be a
condition.  By extending $s$ $\omega$ times if necessary ($R$ is
$\sigma$-closed), we may assume that there is a sequence
$u=\{(p_n,q_n)\}_{n=0}^{\infty}$ in $\omega_2\times Col$ such that $s$
forces $\dot{u}=u$, such that for every $n$, $p_n\in A_s$,
$(p_n,q_n)\in B_s$, that $p_0>_s p_1>_s\cdots >_s p_n>\cdots$ is a
descending chain in $(A_s,<_s)$ and that $q_0 \subset q_1 \subset \dots
\subset q_n \subset \dots.$

Let $p$ be an ordinal greater than $\sup A_s$, let $q\supseteq
\bigcup^\infty_{n=0}q_n$ be such that
$p\in \range(q)\subseteq A_s\cup\{p\}$, let $A_r=A_s\cup\{p\}$,
$B_r=B_s\cup\{(p,q)\}$, $C_r=C_s$, and let  $<_r$ be the partial order
of $A_r$ that extends $<_s$ by making $p$ a lower bound of
$\{p_n\}_{n=0}^{\infty}$.
Finally, let $F_r(\vec a,(p,q))=0$ for all $\vec a \in C_r$
and $r=((A_r,<_r),B_r,C_r,F_r).$

We shall show that for every $\vec a\in C_s,$ $p$ is not a lower bound
of $\vec a.$ This implies that $p \perp_r \vec a$ and (5) and (6) follow.
Since making $p$ a lower bound of $\{p_n\}_n$ does not effect the relation
$\perp$ on $A_s,$ we'll have (8) and hence $r<_R s.$ In $r,$ $(p,q)$ is
a lower bound of $u.$

Thus let $\vec a=\{a_k\}_k \in C_s.$ We claim that
$$\exists k\, \forall n\, a_k \nless_s p_n.$$
This implies that $a_k \nless_r p$ and hence $p$ is not a lower bound
of $\vec a.$

If $p_n < a_0$ for all $n$ then we let $k=0$ because then $a_0 \nless_s p_n$
for all $n.$

Otherwise let $N$ be the least $N$ such that $p_N \ge a_0,$ and let
$m=F_s(\vec a,(p_N,q_N)).$ Either $a_m \nless_s p_n$ for all $n$ and we
are done (with $k=m)$ or else $a_m <_s p_M$ for some $M\ge N.$ By (*) there 
exists some $k$ such that $p_M \perp_s a_k$ and hence $a_k \nless_s p_n$
for all $n.$

\enddemo

Finally, we shall prove that $B(P)$ is not countably closed.

\proclaim
{Lemma 6}  The complete Boolean algebra $B(P)$ does not have a dense
$\sigma$-closed subset.
\endproclaim

\demo
{Proof}  Assume that $B(P)$ does have a dense $\sigma$-closed subset $D$.
For $a,b\in P$, we define
$$
	a\prec b\quad\text{ if }\quad a<_P b\quad\text{ and }\quad
	\exists\; d\in D\quad
	\text{ such that }\quad a<_B d<_B b.
$$
The relation $\prec$ is a partial ordering of $P$, $(P,\prec)$ is
$\sigma$-closed, $a\prec b$ implies $a<_P b$ and for every $a\in P$
there is some $b\in P$ such that $b\prec a$.  

Toward a contradiction,
let $s$ be a condition and assume that $s$ forces 
the preceding statement. For each $\alpha< \omega_2$, there exist a condition 
$s_{\alpha}$ stronger than $s$, and a descending chain 
$\{c^\alpha_n\}_n$ in $A_{s_\alpha}$
such that $c^\alpha_0 \ge \alpha$ and that
for every $n,$ $s_{\alpha}\Vdash c^{\alpha}_{n+1}\prec c^{\alpha}_n$.

By a $\Delta$-system argument we find among these a countable sequence
$r_n=s_{\alpha_n}=((A_n,<_n),B_n,C_n,F_n)$ and a set $D$ such that for
every $m$ and $n$ wih $m<n$ we have $D=A_m \cap A_n,$ $\sup D<\min(A_m-D)$,
$\sup A_m<\min(A_n-D)$, $<_m$ and $<_n$ agree on $D$, 
$\perp_m$ and $\perp_n$ agree on $D$, $B_m\cap(D\times~
Col(D))=B_n\cap(D\times Col(D))$, $C_m\cap D^\omega=C_n\cap D^\omega$,
and $F_m(\vec a,(p,q))=F_n(\vec a,(p,q))$ whenever $\vec a\in C_m\cap
D^\omega$ and $(p,q)\in B_m\cap(D\times Col(D)).$
Moreover, there exists a mappings $\pi_{mn}$ of $A_m$ onto $A_n$ that is an
isomorphism between $(r_m,\{c^{\alpha_m}_k\}_k)$ and 
$(r_n,\{c^{\alpha_n}_k\}_k)$ and is the identity on $D.$
We also let $\pi_{nm}=\pi_{mn}\!^{-1}, \pi_{mm}=id$ and assume that the
$\pi_{mn}$ form a commutative system. Note that for every $n$ and $k$,
$c^{\alpha_n}_k\notin D.$

For each $n$ and $k$, let $a^n_k=c^{\alpha_n}_{2k}$ and $b^n_k=
c^{\alpha_n}_{2k+1}.$ Let $\vec u= \{u_n\}_n$ be the ``diagonal
sequence''
$$
u_{2n}=a^n_n,\quad u_{2n+1}=b^n_n.
$$
\enddemo

We shall find a condition $t=((A_t,<_t),B_t,C_t, F_t)$ stronger than 
all $r_n$ such that the diagonal sequence $\vec u$ is a descending
chain and belongs to $C_t.$ Since $t \Vdash b^n_n \prec a^n_n$ for
every $n,$ it forces that $(P,\prec)$ is not $\sigma$-closed. This will
complete the proof.

To construct $t$ we first let $A_t=\bigcup^\infty_{n=0}A_n$ and 
$B_t=\bigcup^\infty_{n=0}B_n.$ Let
$<_t$ be the minimal partial ordering extending 
$\bigcup^\infty_{n=0}<_n$ such that for every $n,$ $a^{n+1}_{n+1}<_t b^n_n.$
Before proceeding to define $C_t$ and $F_t$ we shall prove some
properties of $(A_t,<_t).$

\proclaim
{Lemma 7} (i) Let $m<n$ and let $y\in A_m-D$ and $x\in A_n-D.$ If $x<_ty$
then $x\le_n a^n_n$ and $b^m_m\le_m y.$ If $x$ and $y$ are compatible in
$<_t$ then $b^m_m\le_m y.$

(ii) For all $m$ and $n,$ if $x\in A_n$ and $y\in A_m$ and if $x<_t y$ then
$x<_n \pi_{mn}y$ (and $\pi_{nm}x<_m y$). In particular, if $x,y\in A_n$ then
$x<_t y$ if and only if $x<_n y.$

(iii) For all $m$ and $n,$ if $x\in A_n$ and $y\in A_m$ and if $x$ and $y$
are compatible in $<_t$ then $x$ and $\pi_{mn}y$ are compatible in $<_n$
(and $\pi_{nm}x$ and $y$ are compatible in $<_m$). In particular, if
$x,y\in A_n$ then $x\perp_t y$ if and only if $x\perp_n y.$
\endproclaim

\demo {Proof} (i) The first statement is an obvious consequence of the
definition of $<_t,$ and the second follows because any $z$ such that
$z\le_t x$ is in some $A_k -D$ where $k\ge n.$

(ii) Let $x\in A_n$ and $y\in A_m$ and let $x<_t y.$ First assume that 
$y\notin D$ (and so $x\notin D.)$ Necessarily, $m\le n$ and if $m=n$ 
then clearly $x<_ny.$ Thus consider $m<n.$ By (i) $x\le_n a^n_n <_n
b^n_m=\pi_{mn}(b^m_m)\le_n\pi_{mn}y.$

Now assume that $y\in D$ and proceed by induction on $x.$ If $x\in D$
then $x<_n y.$ If $x\notin D$ then either $x<_n y$ or there exists 
some $z\notin D$ such that $x<_t z<_t y,$ and by the induction hypothesis
$z<_k \pi_{mk}y$ (where $z\in A_k).$ Applying the preceding paragraph to
$x$ and $z$ we get $\pi_{nk}x<_k z$ and hence $\pi_{nk}x<_k \pi_{mk}y.$
The statement now follows.

(iii) Let $x\in A_n$ and $y\in A_m$ and let $z\in A_k$ be such that
$z<_t x$ and $z<_t y.$ By (ii) we have $\pi_{kn}z<_n x$ and $\pi_{km}z<_m y.$
Hence $\pi_{kn}z=\pi_{mn}\pi_{km}z<_n \pi_{mn}y.$ The second statement
follows from this and from the second statement of (ii). 

\enddemo

Lemma 7 guarantees that $t$ will be stronger than every $r_n.$ Another
consequence is that if $\vec a\in C_n$ then $\vec a$ has no lower bound
in $<_t$: if $x\in A_m$ were a lower bound then $\pi_{mn}x$ would be
a lower bound in $<_n.$

Let $C_t=\bigcup^\infty_{n=0}C_n\cup \{\vec u\}.$ Every sequence in $C_t$
is a descending chain in $<_t$ without a lower bound (clearly, $\vec u$ has
no lower bound).

\proclaim
{Lemma 8} For all $k$ and $n,$ if $(p,q)\in B_k -B_n$ and if $(p',q')\in B_t$
is stronger than $(p,q)$ then $(p',q')\in B_k-B_n.$
\endproclaim

\demo {Proof} Since $(p,q)\notin B_n,$ we have either $\range(q)\not\subseteq D$
or $p\notin D,$ in which case $p\in\range(q)$ by (4) and again $\range(q)\not
\subseteq D.$ Since $q\subseteq q'$ it must be the case 
that $(p',q')\in B_k-B_n$.

We shall now define $F_t$ so that $F_t\supset\bigcup^\infty_{n=0}F_n$ and 
verify (6). This will complete the proof.

First we let $F_t(\vec a,(p,q))=F_n(\vec a,(p,q))$ whenever the right-hand
side is defined; we have to show that (6) holds in $t.$ Let $m=F_n(\vec a,
(p,q))$ and let $(p',q')\in B_k$ be stronger than $(p,q).$ It follows from
Lemma 8 that $(p,q)\in B_k.$ Now $(\pi_{kn}p',\pi_{kn}q')$ is stronger than 
$(\pi_{kn}p,\pi_{kn}q)=(p,q)$ and (*) holds for $\pi_{kn}p'$ in $r_n.$
If $p'<_t a_m$ then by Lemma 7 $\pi_{kn}p'<_n a_m$ and hence $\pi_{kn}p'
\perp_n \vec a.$ By Lemma 7 again, $p'\perp_t\vec a.$

Next, let $\vec a$ and $(p,q)$ be such that $\vec a\in C_n - C_k,$ $(p,q)\in
B_k - B_n$ and $p\ge a_0.$ If $k<n,$ we have $\pi_{kn}p\ge p\ge a_0$ and we
let $F_t(\vec a,(p,q))=F_n(\vec a,(\pi_{kn}p,\pi_{kn}q)).$ To verify (6), let
$m=F_t(\vec a,(p,q))$ and let $(p',q')\in B_t$ be stronger than $(p,q).$ 
By Lemma 8 $(p',q')\in B_k,$ and $(\pi_{kn}p',\pi_{kn}q')$ is stronger 
(in $r_n$) than $(\pi_{kn}p,\pi_{kn}q).$ If $p'<_t a_m$ then by Lemma 7 
$\pi_{kn}p'<_n a_m$ and so $\pi_{kn}p'\perp_n \vec a.$ By Lemma 7 again, 
$p'\perp_t \vec a.$

If $k>n,$ we let $F_t(\vec a,(p,q))$ be the least $m$ such that $a_m\notin D$
and that $b^n_n \not\le a_m.$ To verify (6), let $(p',q')\in B_t$ be stronger 
than $(p,q).$ If $p'\in D$ then $p'\nless_t a_m$ and if $p'\notin D$ then
by Lemma 7(i) $p'\perp_t a_m.$ In either case, (6) is satisfied.

Finally, we define $F_t(\vec u,(p,q)).$ Thus let $(p,q)\in B_t$ be such that
$p\ge u_0.$ Since $u_0=a^0_0\notin D,$ we have $p\notin D$. Let $n$ be the $n$
such that $p\in A_n.$ We let $F_t(\vec u,(p,q))=2n+2.$ That is, the chosen
$u_m$ is $u_{2n+2}=a^{n+1}_{n+1}.$ To verify (6), let $(p',q')\in B_t$ be
stronger than $(p,q).$ Since $p\in A_n -D,$ by Lemma 8 we have 
$(p',q')\in B_n$ and therefore $p'\in A_n -D.$ But $a^{n+1}_{n+1}\in A_{n+1}-D$
and so $p'\nless_t a^{n+1}_{n+1}.$ Therefore (6) holds.
\enddemo

\enddocument